\documentclass[journal]{IEEEtran}

\usepackage[utf8]{inputenc}
\usepackage{color}
\usepackage{array}
\usepackage{verbatim}
\usepackage{float}
\usepackage{amsmath}
\usepackage{amsthm}
\usepackage{amssymb}
\usepackage{graphicx}
\usepackage{longtable}
\usepackage[unicode=true,
 bookmarks=false,
 breaklinks=false,pdfborder={0 0 1},colorlinks=false]
 {hyperref}
\hypersetup{
 colorlinks,bookmarksopen,bookmarksnumbered,citecolor=blue,urlcolor=blue}
\usepackage{cite}
\makeatletter
 \let\oldforeign@language\foreign@language
 \DeclareRobustCommand{\foreign@language}[1]{%
   \lowercase{\oldforeign@language{#1}}}

 \let\oldforeign@language\foreign@language
 \DeclareRobustCommand{\foreign@language}[1]{%
   \lowercase{\oldforeign@language{#1}}}

\ifCLASSINFOpdf
\else
\fi

\makeatletter

\def\ps@IEEEtitlepagestyle{%
	\def\@oddhead{\parbox[t][\height][t]{\textwidth}{\centering
			\scriptsize Personal use of this material is permitted. Permission from the author(s) and/or copyright holder(s), must be obtained for all other uses. Please contact us and provide details if you believe this document breaches copyrights.\\
			\noindent\makebox[\linewidth]{}
		}\hfil\hbox{}}%
	\def\@evenhead{\scriptsize\thepage \hfil \leftmark\mbox{}}%
	\def\@oddfoot{\parbox[t][\height][t]{\textwidth}{
			\vspace{-20pt}{\rule{\textwidth}{0.4pt}}\\\footnotesize \underline{To cite this article:}
			{\bf{\textcolor{red}{H. A. Hashim, L. J. Brown, and K. McIsaac, "Guaranteed Performance of Nonlinear Pose Filter on SE(3),"  in Proceedings of the 2019 American Control Conference (ACC), Philadelphia, PA, USA, 2019, pp. 1108-1113.}}} doi: \href{https://doi.org/10.23919/ACC.2019.8814878}{10.23919/ACC.2019.8814878} 
			\noindent\makebox[\linewidth]
		}\hfil\hbox{}}%
	\def\@evenfoot{\MYfooter}}

\makeatother
\newtheorem{lem}{Lemma}

\newtheorem{prop}{Proposition}
\newtheorem{thm}{Theorem}
\newtheorem{rem}{Remark}

\newtheorem{assum}{Assumption}

\begin{document}
\bstctlcite{IEEEexample:BSTcontrol}
%
%
%
%
%
%
%

\title{Guaranteed Performance of Nonlinear Pose Filter on SE(3)}

\author{Hashim~A.~Hashim, Lyndon J. Brown, and~Kenneth McIsaac
\thanks{H. A. Hashim, L. J. Brown and K. McIsaac are with the Department of Electrical and Computer Engineering,
University of Western Ontario, London, ON, Canada, N6A-5B9, e-mail: hmoham33@uwo.ca, lbrown@uwo.ca and kmcisaac@uwo.ca.}
}


\markboth{--,~Vol.~-, No.~-, \today}{Hashim \MakeLowercase{\textit{et al.}}: Guaranteed Performance of Nonlinear Pose Filter on SE(3)}
\markboth{}{Hashim \MakeLowercase{\textit{et al.}}: Guaranteed Performance of Nonlinear Pose Filter on SE(3)p}

\maketitle

\begin{abstract}
This paper presents a novel nonlinear pose filter evolved directly
on the Special Euclidean Group $\mathbb{SE}\left(3\right)$ with guaranteed
characteristics of transient and steady-state performance. The above-mention
characteristics can be achieved by trapping the position error and
the error of the normalized Euclidean distance of the attitude in
a given large set and guiding them to converge systematically to a
small given set. The error vector is proven to approach the origin
asymptotically from almost any initial condition. The proposed filter
is able to provide a reliable pose estimate with remarkable convergence
properties such that it can be fitted with measurements obtained from
low-cost measurement units. Simulation results demonstrate high convergence
capabilities and robustness considering large error in initialization
and high level of uncertainties in measurements.  
\end{abstract}

%

\IEEEpeerreviewmaketitle{}

\section{Introduction}

Pose of a rigid-body in 3D space can be described by two components:
orientation and translation. A reasonable pose estimation of the rigid-body
in 3D space is crucial for robotics and engineering applications,
such as space crafts, unmanned aerial and underwater vehicles, satellites,
etc. The orientation (attitude) can be established using statical
methods, such as QUEST \cite{shuster1981three} and singular value
decomposition (SVD) \cite{markley1988attitude}, which utilize a set
of known vectors in the inertial-frame and their measurements in the
body-frame. However, body-frame measurements are contaminated with
noise and bias components \cite{hashim2018SO3Stochastic,hashim2018Conf1,hashim2019SO3Det}
causing the static estimation algorithms in \cite{shuster1981three,markley1988attitude}
to produce unsatisfactory results.

The attitude can be estimated through Gaussian filters which often
consider unit-quaternion in the representation, such as Kalman filter
(KF) \cite{choukroun2006novel}, extended KF (EKF) \cite{lefferts1982kalman},
and multiplicative EKF (MEKF) \cite{markley2003attitude}. However,
to successfully address the nonlinear nature of the attitude problem
a nonlinear deterministic filter evolved directly on the Special Orthogonal
Group $\mathbb{SO}\left(3\right)$ can be used \cite{hashim2018SO3Stochastic,hashim2018Conf1,hashim2019SO3Det,mahony2008nonlinear,liu2018complementary,wu2015globally,bohn2014almost}.
As a matter of fact, nonlinear deterministic attitude filters are
simpler in derivation, require less computational power, and demonstrate
better tracking performance in comparison with Gaussian filters \cite{hashim2018SO3Stochastic}.
It should be remarked that attitude is a major part of the pose problem.
As such, the pose filtering problem is better addressed in the nonlinear
sense.

The pose filter could be developed based on the measurements obtained
from inertial measurement units (IMUs) along with landmark measurements
collected by a vision system. The observer in \cite{baldwin2007complementary}
was evolved directly on the Special Euclidean Group $\mathbb{SE}\left(3\right)$
and, while it required pose reconstruction, it was subsequently adjusted
in \cite{baldwin2009nonlinear,hua2011observer} to function based
solely on a set of vectorial measurements. A recent nonlinear stochastic
pose observer on $\mathbb{SE}\left(3\right)$ applicable for measurements
obtained from low-cost measurement units is proposed \cite{hashim2018SE3Stochastic}.
Although the filters discussed in \cite{baldwin2007complementary,baldwin2009nonlinear,hua2011observer}
are simple in design, they are highly sensitive to the uncertain measurements.
Moreover, there is no guarantee that the tracking error will behave
according to the predefined dynamic constraints of the transient and
steady-state performance. Prescribed performance can be defined as
a process of systematic convergence of the error from a large known
set to a small known set guided by the prescribed performance function
(PPF) \cite{bechlioulis2008robust}. The constrained error is transformed
to unconstrained form termed transformed error. The remarkable advantage
offered by PPF could be utilized in control and filtering design process
of two degree of freedom planar robot \cite{bechlioulis2008robust},
uncertain multi-agent system \cite{hashim2017neuro} and other applications. 

This paper presents a robust nonlinear pose filter on $\mathbb{SE}\left(3\right)$
that satisfies predefined characteristics of transient and steady-state
measures. The error initially starts within a predefined large set
and is forced to decrease systematically to a given small set with
the aid of the transformed error. The error of the homogeneous transformation
matrix asymptotically approaches the identity, as the transformed
error approaches the origin and vice versa. The filter is guaranteed
to demonstrate fast convergence and robustness against high level
of uncertainties in the measurements from almost any initial condition. 

The rest of the paper is organized as follows: Section \ref{sec:SE3PPF_Math-Notations}
provides the preliminaries of $\mathbb{SO}\left(3\right)$ and $\mathbb{SE}\left(3\right)$.
Vector measurements are presented and the pose problem is formulated
in terms of prescribed performance in Section \ref{sec:SE3PPF_Problem-Formulation-in}.
The nonlinear pose filter and the stability analysis are laid out
in Section \ref{sec:SE3PPF-Filters}. Section \ref{sec:SO3PPF_Simulations}
illustrates the fast convergence and robustness of the proposed filter.
Finally, Section \ref{sec:SO3PPF_Conclusion} concludes the work.

\section{Preliminaries of $\mathbb{SE}\left(3\right)$ \label{sec:SE3PPF_Math-Notations}}

In this paper $\mathbb{R}_{+}$, $\mathbb{R}^{n}$ and $\mathbb{R}^{n\times m}$
denote the set of non-negative real numbers, real $n$-dimensional
space column vector, and real $n\times m$ dimensional space, respectively.
The Euclidean norm of $x\in\mathbb{R}^{n}$ is defined by $\left\Vert x\right\Vert =\sqrt{x^{\top}x}$.
$\mathbf{I}_{n}$ refers to an $n$-by-$n$ identity matrix, where
$\underline{\mathbf{0}}_{n}$ is a zero column vector. Define $\mathbb{SO}\left(3\right)$
as the Special Orthogonal Group. The orientation of a rigid-body in
space, also known as attitude matrix $R$, is given by
\[
\mathbb{SO}\left(3\right)=\left\{ \left.R\in\mathbb{R}^{3\times3}\right|RR^{\top}=R^{\top}R=\mathbf{I}_{3}\text{, }{\rm det}\left(R\right)=+1\right\} 
\]
where $\mathbf{I}_{3}$ is a $3$-by-$3$ identity matrix and ${\rm det\left(\cdot\right)}$
is a determinant of the matrix. Define $\mathbb{SE}\left(3\right)$
as the Special Euclidean Group with $\mathbb{SE}\left(3\right)$ being
defined by
\[
\mathbb{SE}\left(3\right)=\left\{ \left.\boldsymbol{T}\in\mathbb{R}^{4\times4}\right|R\in\mathbb{SO}\left(3\right),P\in\mathbb{R}^{3}\right\} 
\]
with $\boldsymbol{T}\in\mathbb{SE}\left(3\right)$ being the homogeneous
transformation matrix that describes the pose of the rigid-body as
follows
\begin{equation}
\boldsymbol{T}=\left[\begin{array}{cc}
R & P\\
\underline{\mathbf{0}}_{3}^{\top} & 1
\end{array}\right]\in\mathbb{SE}\left(3\right)\label{eq:SE3STCH_T_matrix}
\end{equation}
with $P\in\mathbb{R}^{3}$ and $R\in\mathbb{SO}\left(3\right)$ standing
for position and attitude of the rigid-body in space, respectively,
and $\underline{\mathbf{0}}_{3}^{\top}$ being a zero row. The Lie-algebra
related to $\mathbb{SO}\left(3\right)$ is termed $\mathfrak{so}\left(3\right)$
and is given by
\[
\mathfrak{so}\left(3\right)=\left\{ \left.B\in\mathbb{R}^{3\times3}\right|B^{\top}=-B\right\} 
\]
where $B$ is a skew symmetric matrix. Let the map $\left[\cdot\right]_{\times}:\mathbb{R}^{3}\rightarrow\mathfrak{so}\left(3\right)$
be
\[
\left[\beta\right]_{\times}=\left[\begin{array}{ccc}
0 & -\beta_{3} & \beta_{2}\\
\beta_{3} & 0 & -\beta_{1}\\
-\beta_{2} & \beta_{1} & 0
\end{array}\right]\in\mathfrak{so}\left(3\right),\hspace{1em}\beta=\left[\begin{array}{c}
\beta_{1}\\
\beta_{2}\\
\beta_{3}
\end{array}\right]
\]
Define $\left[\beta\right]_{\times}\vartheta=\beta\times\vartheta$
with $\times$ being the cross product for all $\beta,\vartheta\in\mathbb{R}^{3}$.
For any $\mathcal{Y}=\left[y_{1}^{\top},y_{2}^{\top}\right]^{\top}$
with $y_{1},y_{2}\in\mathbb{R}^{3}$, the wedge map $\left[\cdot\right]_{\wedge}:\mathbb{R}^{6}\rightarrow\mathfrak{se}\left(3\right)$
is given by
\[
\left[\mathcal{Y}\right]_{\wedge}=\left[\begin{array}{cc}
\left[y_{1}\right]_{\times} & y_{2}\\
\underline{\mathbf{0}}_{3}^{\top} & 0
\end{array}\right]\in\mathfrak{se}\left(3\right)
\]
Let $\mathfrak{se}\left(3\right)$ be the Lie algebra of $\mathbb{SE}\left(3\right)$
defined by{\small{}
	\begin{align*}
	\mathfrak{se}\left(3\right) & =\left\{ \left.\left[\mathcal{Y}\right]_{\wedge}\in\mathbb{R}^{4\times4}\right|\exists y_{1},y_{2}\in\mathbb{R}^{3}:\left[\mathcal{Y}\right]_{\wedge}=\left[\begin{array}{cc}
	\left[y_{1}\right]_{\times} & y_{2}\\
	\underline{\mathbf{0}}_{3}^{\top} & 0
	\end{array}\right]\right\} 
	\end{align*}
}Consider the inverse map of $\left[\cdot\right]_{\times}$ such that
$\mathbf{vex}:\mathfrak{so}\left(3\right)\rightarrow\mathbb{R}^{3}$
\begin{equation}
\mathbf{vex}\left(\left[\beta\right]_{\times}\right)=\beta\in\mathbb{R}^{3}\label{eq:SE3STCH_VEX}
\end{equation}
The anti-symmetric projection operator on the Lie-algebra $\mathfrak{so}\left(3\right)$
is denoted by $\boldsymbol{\mathcal{P}}_{a}$ with $\boldsymbol{\mathcal{P}}_{a}:\mathbb{R}^{3\times3}\rightarrow\mathfrak{so}\left(3\right)$
such that
\begin{equation}
\boldsymbol{\mathcal{P}}_{a}\left(M\right)=\frac{1}{2}\left(M-M^{\top}\right)\in\mathfrak{so}\left(3\right),\,M\in\mathbb{R}^{3\times3}\label{eq:SE3STCH_Pa}
\end{equation}
The normalized Euclidean distance of $R\in\mathbb{SO}\left(3\right)$
is
\begin{equation}
\left\Vert R\right\Vert _{I}=\frac{1}{4}{\rm Tr}\left\{ \mathbf{I}_{3}-R\right\} \in\left[0,1\right]\label{eq:SE3STCH_Ecul_Dist}
\end{equation}
where ${\rm Tr}\left\{ \cdot\right\} $ is a trace of a matrix. The
following mathematical identity will be used in the filter derivation
\begin{align}
{\rm Tr}\left\{ A\left[\beta\right]_{\times}\right\} = & {\rm Tr}\left\{ \boldsymbol{\mathcal{P}}_{a}\left(A\right)\left[\beta\right]_{\times}\right\} =-2\mathbf{vex}\left(\boldsymbol{\mathcal{P}}_{a}\left(A\right)\right)^{\top}\beta\nonumber \\
& \hspace{1em}A\in\mathbb{R}^{3\times3},\beta\in{\rm \mathbb{R}}^{3}\label{eq:SO3PPF_Identity7}
\end{align}

\section{Problem Formulation \label{sec:SE3PPF_Problem-Formulation-in}}

The aim of this section is to present the pose problem, introduce
pose measurements, and reformulate the problem in terms of prescribed
performance.

\subsection{Pose Dynamics and Measurements\label{subsec:SE3PPF_Pose-Kinematics}}

\begin{figure}[h]
	\centering{}\includegraphics[scale=0.35]{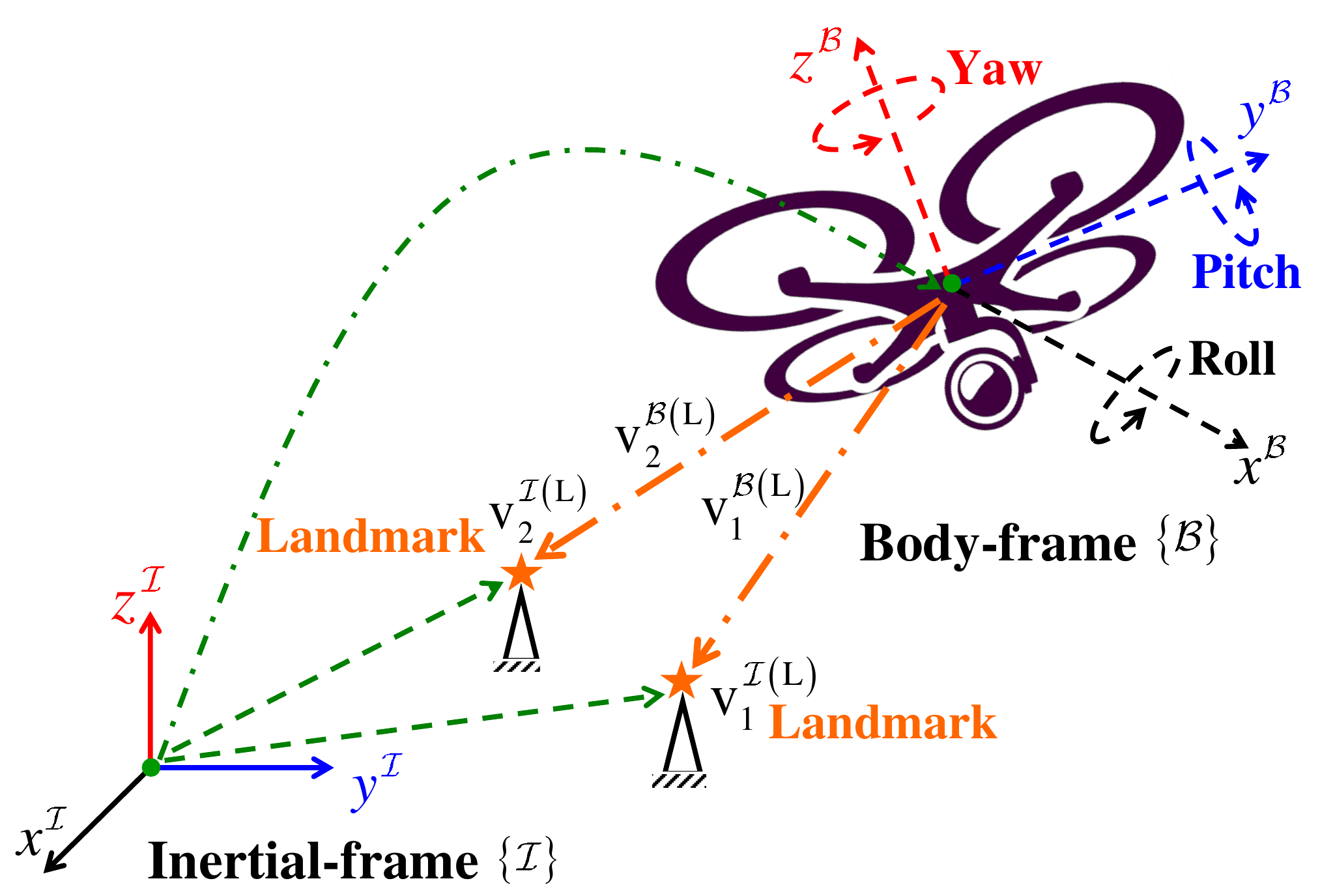}\caption{Pose estimation problem of a rigid-body in 3D space.}
	\label{fig:SE3PPF_1} 
\end{figure}

The pose of a rigid-body is determined by its attitude and position.
The attitude of a rigid-body is given by $R\in\mathbb{SO}\left(3\right)$
with $R\in\left\{ \mathcal{B}\right\} $ while the position is defined
by $P\in\mathbb{R}^{3}$ with $P\in\left\{ \mathcal{I}\right\} $.
The pose estimation of a rigid-body illustrated in Fig. \ref{fig:SE3PPF_1}
can be described by the following homogeneous transformation matrix
\begin{equation}
\boldsymbol{T}=\left[\begin{array}{cc}
R & P\\
\underline{\mathbf{0}}_{3}^{\top} & 1
\end{array}\right]\in\mathbb{SE}\left(3\right)\label{eq:SE3PPF_T_matrix2}
\end{equation}

Define the superscripts $\mathcal{B}$ and $\mathcal{I}$ as the components
associated with the body-frame and inertial-frame, respectively. The
attitude can be expressed through $N_{{\rm R}}$ known measured vectors
in the body-frame and those vectors are known in the inertial frame.
The $j$th vector measurement in the body-frame is defined by
\begin{equation}
{\rm v}_{j}^{\mathcal{B}\left({\rm R}\right)}=R^{\top}{\rm v}_{j}^{\mathcal{I}\left({\rm R}\right)}+{\rm b}_{j}^{\mathcal{B}\left({\rm R}\right)}+\omega_{j}^{\mathcal{B}\left({\rm R}\right)}\in\mathbb{R}^{3}\label{eq:SE3STCH_Vect_R}
\end{equation}
where ${\rm v}_{j}^{\mathcal{I}\left({\rm R}\right)}$ is a known
vector, ${\rm b}_{j}^{\mathcal{B}\left({\rm R}\right)}$ is unknown
bias, and $\omega_{j}^{\mathcal{B}\left({\rm R}\right)}$ is unknown
random noise associated with the $j$th measurement for all $j=1,2,\ldots,N_{{\rm R}}$.
The vectors ${\rm v}_{j}^{\mathcal{I}\left({\rm R}\right)}$ and ${\rm v}_{j}^{\mathcal{B}\left({\rm R}\right)}$
in \eqref{eq:SE3STCH_Vect_R} can be normalized as follows
\begin{equation}
\upsilon_{j}^{\mathcal{I}\left({\rm R}\right)}=\frac{{\rm v}_{j}^{\mathcal{I}\left({\rm R}\right)}}{\left\Vert {\rm v}_{j}^{\mathcal{I}\left({\rm R}\right)}\right\Vert },\hspace{1em}\upsilon_{j}^{\mathcal{B}\left({\rm R}\right)}=\frac{{\rm v}_{j}^{\mathcal{B}\left({\rm R}\right)}}{\left\Vert {\rm v}_{j}^{\mathcal{B}\left({\rm R}\right)}\right\Vert }\label{eq:SE3STCH_Vector_norm}
\end{equation}
The position of a moving body can be reconstructed provided that $R$
is known and there exist $N_{{\rm L}}$ known landmarks. The $j$th
vector measurement in the body-frame is given by
\begin{equation}
{\rm v}_{j}^{\mathcal{B}\left({\rm L}\right)}=R^{\top}\left({\rm v}_{j}^{\mathcal{I}\left({\rm L}\right)}-P\right)+{\rm b}_{j}^{\mathcal{B}\left({\rm L}\right)}+\omega_{j}^{\mathcal{B}\left({\rm L}\right)}\in\mathbb{R}^{3}\label{eq:SE3STCH_Vec_Landmark}
\end{equation}
with ${\rm v}_{j}^{\mathcal{I}\left({\rm L}\right)}$ being a known
feature, ${\rm b}_{j}^{\mathcal{B}\left({\rm L}\right)}$ being unknown
bias, and $\omega_{j}^{\mathcal{B}\left({\rm R}\right)}$ being unknown
random noise of the $j$th measurement for all $j=1,2,\ldots,N_{{\rm L}}$. 
\begin{assum}
	\label{Assum:SE3STCH_1} (Pose observability) At least three non-collinear
	vectors in \eqref{eq:SE3STCH_Vector_norm} and one landmark in \eqref{eq:SE3STCH_Vec_Landmark}
	must be available in order to extract the pose of a rigid-body. If
	$N_{{\rm R}}=2$, the third non-collinear vector can be obtained by
	$\upsilon_{3}^{\mathcal{I}\left({\rm R}\right)}=\upsilon_{1}^{\mathcal{I}\left({\rm R}\right)}\times\upsilon_{2}^{\mathcal{I}\left({\rm R}\right)}$
	and $\upsilon_{3}^{\mathcal{B}\left({\rm R}\right)}=\upsilon_{1}^{\mathcal{B}\left({\rm R}\right)}\times\upsilon_{2}^{\mathcal{B}\left({\rm R}\right)}$. 
\end{assum}
For simplicity, ${\rm v}_{j}^{\mathcal{B}\left({\rm R}\right)}$ and
${\rm v}_{j}^{\mathcal{B}\left({\rm L}\right)}$ are assumed to be
free of noise and bias components in the stability analysis. In the
simulation section, however, noise and bias present in the measurements
are taken into consideration. The pose dynamics of a rigid-body are
defined by
\begin{align}
\dot{\boldsymbol{T}} & =\boldsymbol{T}\left[\mathcal{Y}\right]_{\wedge}\label{eq:SE3PPF_T_Dynamics}
\end{align}
where $\dot{R}=R\left[\Omega\right]_{\times}$, $\dot{P}=RV$, $\mathcal{Y}=\left[\Omega^{\top},V^{\top}\right]^{\top}\in\mathbb{R}^{6}$
is a group velocity vector with $\Omega\in\mathbb{R}^{3}$ and $V\in\mathbb{R}^{3}$
being the true angular and translational velocities, respectively.
The measured velocity vector is defined by
\begin{align}
\mathcal{Y}_{m} & =\mathcal{Y}+b+\omega\in\left\{ \mathcal{B}\right\} \label{eq:SE3PPF_Angular}
\end{align}
where $\mathcal{Y}_{m}=\left[\Omega_{m}^{\top},V_{m}^{\top}\right]^{\top}$,
$b=\left[b_{\Omega}^{\top},b_{V}^{\top}\right]^{\top}$, and $\omega=\left[\omega_{\Omega}^{\top},\omega_{V}^{\top}\right]^{\top}$,
with $b_{\Omega},b_{V}\in\mathbb{R}^{3}$ being unknown constant bias
and $\omega_{\Omega},\omega_{V}\in\mathbb{R}^{3}$ being unknown random
noise attached to the measurements. In this section, in the interest
of simplicity it is assumed that $\omega=\underline{\mathbf{0}}_{6}$,
while in the implementation $\omega\neq\underline{\mathbf{0}}_{6}$.
From the identity in \eqref{eq:SO3PPF_Identity7}, the dynamics of
the normalized Euclidean distance are defined by
\begin{align}
||\dot{R}||_{I} & =\frac{1}{2}\mathbf{vex}\left(\boldsymbol{\mathcal{P}}_{a}\left(R\right)\right)^{\top}\Omega\label{eq:SE3PPF_NormR_dynam}
\end{align}
Consequently, the pose kinematics in \eqref{eq:SE3PPF_T_Dynamics}
can be expressed in vector form by{\small{}
	\begin{equation}
	\left[\begin{array}{c}
	||\dot{R}||_{I}\\
	\dot{P}
	\end{array}\right]=\left[\begin{array}{cc}
	\frac{1}{2}\mathbf{vex}\left(\boldsymbol{\mathcal{P}}_{a}\left(R\right)\right)^{\top} & \underline{\mathbf{0}}_{3}^{\top}\\
	\mathbf{0}_{3\times3} & R
	\end{array}\right]\left[\begin{array}{c}
	\Omega_{m}-b_{\Omega}\\
	V_{m}-b_{V}
	\end{array}\right]\label{eq:SE3PPF_T_VEC_Dyn}
	\end{equation}
}{\small\par}

\noindent Define the estimate of the homogeneous transformation matrix
($\hat{\boldsymbol{T}}$) in \eqref{eq:SE3PPF_T_matrix2} by
\begin{equation}
\hat{\boldsymbol{T}}=\left[\begin{array}{cc}
\hat{R} & \hat{P}\\
\underline{\mathbf{0}}_{3}^{\top} & 1
\end{array}\right]\label{eq:SE3PPF_Test_matrix}
\end{equation}
where $\hat{R}$ and $\hat{P}$ are the estimates of $R$ and $P$,
respectively. Define the homogeneous transformation matrix error by
\begin{align}
\tilde{\boldsymbol{T}} & =\hat{\boldsymbol{T}}\boldsymbol{T}^{-1}=\left[\begin{array}{cc}
\tilde{R} & \tilde{P}\\
\underline{\mathbf{0}}_{3}^{\top} & 1
\end{array}\right]\label{eq:SE3PPF_Terr_matrix}
\end{align}

\noindent where $\tilde{R}=\hat{R}R^{\top}$ and $\tilde{P}=\hat{P}-\tilde{R}P$
are the errors in attitude and position, respectively. The objective
of this work is to drive $\hat{\boldsymbol{T}}\rightarrow\boldsymbol{T}$
which ensures that $\tilde{P}\rightarrow\underline{\mathbf{0}}_{3}$,
$\tilde{R}\rightarrow\mathbf{I}_{3}$, and $\tilde{\boldsymbol{T}}\rightarrow\mathbf{I}_{4}$.
The following Lemma \ref{Lemm:SE3PPF_1} is important in the filter
derivation. 
\begin{lem}
	\label{Lemm:SE3PPF_1}Let $R\in\mathbb{SO}\left(3\right)$. Then,
	the following holds:
	\begin{equation}
	||\mathbf{vex}\left(\boldsymbol{\mathcal{P}}_{a}\left(R\right)\right)||^{2}=4\left(1-||R||_{I}\right)||R||_{I}\label{eq:SE3PPF_lemm1_1}
	\end{equation}
	\textbf{Proof. See \nameref{sec:SO3STCH_EXPL_AppendixA}.}
\end{lem}

\subsection{Prescribed Performance \label{subsec:SE3PPF_Prescribed-Performance}}

\begin{figure}[h!]
	\centering{}\includegraphics[scale=0.27]{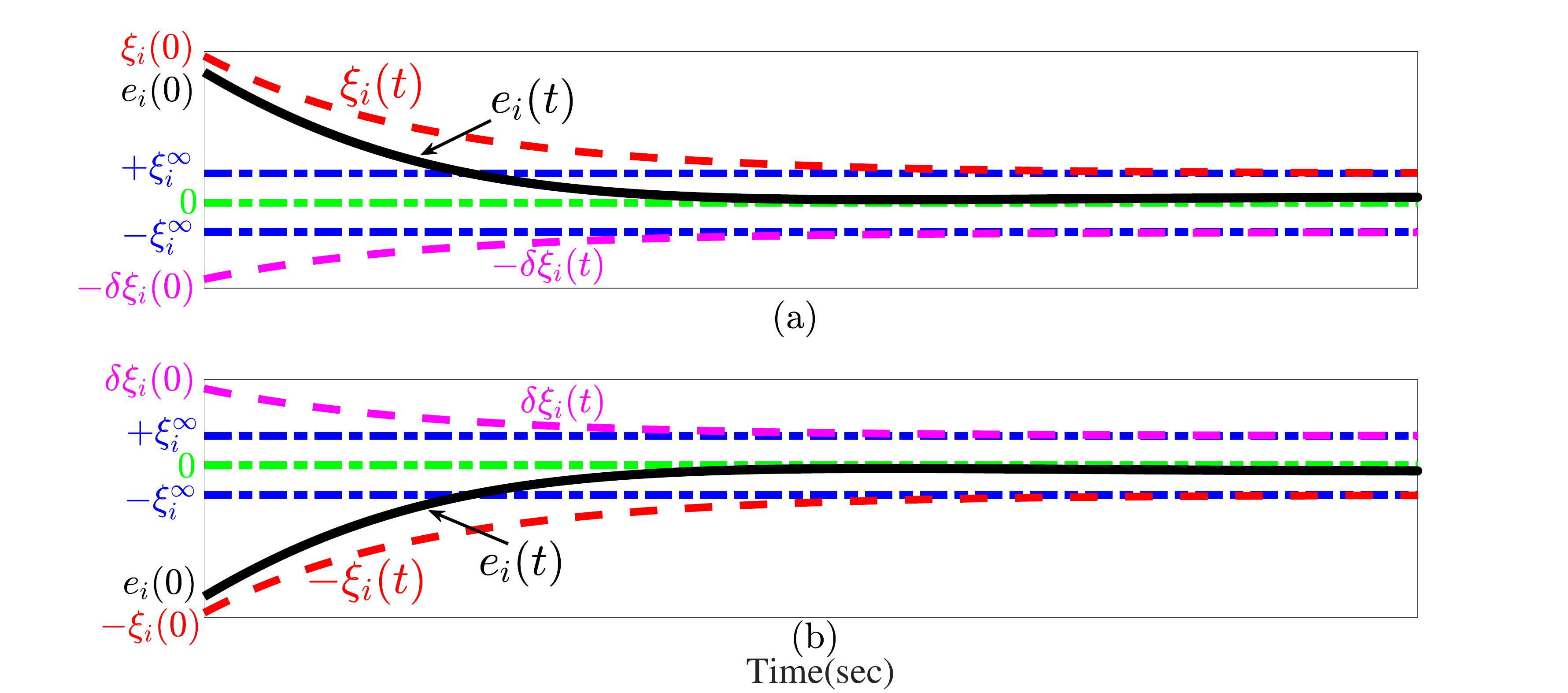} \caption{Systematic convergence of $\boldsymbol{e}_{i}$ with PPF satisfying
		(a) Eq. \eqref{eq:SE3PPF_ePos}; (b) Eq. \eqref{eq:SE3PPF_eNeg}.}
	\label{fig:SO3PPF_2} 
\end{figure}

Let the error in the homogeneous transformation matrix be as in \eqref{eq:SE3PPF_Terr_matrix}.
In view of \eqref{eq:SE3PPF_T_VEC_Dyn}, define the error in vector
form as
\begin{equation}
\boldsymbol{e}=\left[\boldsymbol{e}_{1},\boldsymbol{e}_{2},\boldsymbol{e}_{3},\boldsymbol{e}_{4}\right]^{\top}=\left[||\tilde{R}||_{I},\tilde{P}^{\top}\right]^{\top}\in\mathbb{R}^{4}\label{eq:SE3PPF_Vec_error}
\end{equation}
The aim is to initiate the error within a given large set and reduce
it systematically and smoothly to a given small set using the prescribed
performance function (PPF) \cite{bechlioulis2008robust}. Define the
following PPF \cite{bechlioulis2008robust}
\begin{equation}
\xi_{i}\left(t\right)=\left(\xi_{i}^{0}-\xi_{i}^{\infty}\right)\exp\left(-\ell_{i}t\right)+\xi_{i}^{\infty}\label{eq:SE3PPF_Presc}
\end{equation}
with $\xi_{i}\left(t\right)$ being a time-decreasing positive smooth
function that satisfies $\xi_{i}:\mathbb{R}_{+}\to\mathbb{R}_{+}$.
Also, $\lim\limits _{t\to\infty}\xi_{i}\left(t\right)=\xi_{i}^{\infty}>0$
with $\xi_{i}\left(0\right)=\xi_{i}^{0}$ being the initial value
and the upper bound of $\xi_{i}\left(t\right)$, $\xi_{i}^{\infty}$
being the upper bound of the small set, and the positive constant
$\ell_{i}$ controlling the convergence rate of $\xi_{i}\left(t\right)$
for all $i=1,2,\ldots,4$. Meeting the following conditions is sufficient
to ensure the systematic convergence of $\boldsymbol{e}_{i}\left(t\right)$
within the PPF:
\begin{align}
-\delta\xi_{i}\left(t\right)<\boldsymbol{e}_{i}\left(t\right)<\xi_{i}\left(t\right), & \text{ if }\boldsymbol{e}_{i}\left(0\right)\geq0\label{eq:SE3PPF_ePos}\\
-\xi_{i}\left(t\right)<\boldsymbol{e}_{i}\left(t\right)<\delta\xi_{i}\left(t\right), & \text{ if }\boldsymbol{e}_{i}\left(0\right)<0\label{eq:SE3PPF_eNeg}
\end{align}
such that $1\geq\delta\geq0$. For clarity, let $\xi=\left[\xi_{1},\ldots,\xi_{4}\right]^{\top}$,
$\ell=\left[\ell_{1},\ldots,\ell_{4}\right]^{\top}$, $\xi^{0}=\left[\xi_{1}^{0},\ldots,\xi_{4}^{0}\right]^{\top}$,
and $\xi^{\infty}=\left[\xi_{1}^{\infty},\ldots,\xi_{4}^{\infty}\right]^{\top}$
with $\boldsymbol{e}_{i}:=\boldsymbol{e}_{i}\left(t\right)$ and $\xi_{i}:=\xi_{i}\left(t\right)$
for all $\xi,\ell,\xi^{0},\xi^{\infty}\in\mathbb{R}^{4}$. The systematic
convergence of $\boldsymbol{e}_{i}$ from a known large set to a known
small set is depicted in Fig. \ref{fig:SO3PPF_2}.

\begin{rem}
	\label{SE3PPF_rem3}\cite{hashim2017neuro,hashim2017adaptive,bechlioulis2008robust}
	Knowing the upper bound and the sign of $\boldsymbol{e}_{i}\left(0\right)$
	is sufficient to ensure that the error follows the PPF and that the
	tracking error is maintained within known dynamically reducing boundaries
	for all $t>0$ as depicted in Fig. \ref{fig:SO3PPF_2}. 
\end{rem}
Define the error $\boldsymbol{e}_{i}$ by
\begin{equation}
\boldsymbol{e}_{i}=\xi_{i}\mathcal{Z}\left(\mathcal{E}_{i}\right)\label{eq:SE3PPF_e_Trans}
\end{equation}
where $\mathcal{E}_{i}\in\mathbb{R}$ is an unconstrained transformed
error, and $\mathcal{Z}\left(\mathcal{E}_{i}\right)$ possessing the
properties listed below: 
\begin{enumerate}
	\item[(i)] $\mathcal{Z}\left(\mathcal{E}_{i}\right)$ is a smooth and increasing
	function. 
	\item[(ii)] $\mathcal{Z}\left(\mathcal{E}_{i}\right)$ is constrained such that
	\\
	$-\underline{\delta}_{i}<\mathcal{Z}\left(\mathcal{E}_{i}\right)<\bar{\delta}_{i},{\rm \text{ if }}\boldsymbol{e}_{i}\left(0\right)\geq0$\\
	$-\bar{\delta}_{i}<\mathcal{Z}\left(\mathcal{E}_{i}\right)<\underline{\delta}_{i},{\rm \text{ if }}\boldsymbol{e}_{i}\left(0\right)<0$
	\\
	with $\bar{\delta}_{i},\underline{\delta}_{i}>0$ and $\underline{\delta}_{i}\leq\bar{\delta}_{i}$. 
	\item[(iii)] $\left.\begin{array}{c}
	\underset{\mathcal{E}_{i}\rightarrow-\infty}{\lim}\mathcal{Z}\left(\mathcal{E}_{i}\right)=-\underline{\delta}_{i}\\
	\underset{\mathcal{E}_{i}\rightarrow+\infty}{\lim}\mathcal{Z}\left(\mathcal{E}_{i}\right)=\bar{\delta}_{i}
	\end{array}\right\} {\rm \text{ if }}\boldsymbol{e}_{i}\geq0$ 
	\item[] $\left.\begin{array}{c}
	\underset{\mathcal{E}_{i}\rightarrow-\infty}{\lim}\mathcal{Z}\left(\mathcal{E}_{i}\right)=-\bar{\delta}_{i}\\
	\underset{\mathcal{E}_{i}\rightarrow+\infty}{\lim}\mathcal{Z}\left(\mathcal{E}_{i}\right)=\underline{\delta}_{i}
	\end{array}\right\} {\rm \text{ if }}\boldsymbol{e}_{i}<0$ \\
	such that
\end{enumerate}
\begin{equation}
\mathcal{Z}\left(\mathcal{E}_{i}\right)=\begin{cases}
\frac{\bar{\delta}_{i}\exp(\mathcal{E}_{i})-\underline{\delta}_{i}\exp(-\mathcal{E}_{i})}{\exp(\mathcal{E}_{i})+\exp(-\mathcal{E}_{i})}, & \bar{\delta}_{i}\geq\underline{\delta}_{i}\text{ if }\boldsymbol{e}_{i}\geq0\\
\frac{\bar{\delta}_{i}\exp(\mathcal{E}_{i})-\underline{\delta}_{i}\exp(-\mathcal{E}_{i})}{\exp(\mathcal{E}_{i})+\exp(-\mathcal{E}_{i})}, & \underline{\delta}_{i}\geq\bar{\delta}_{i}\text{ if }\boldsymbol{e}_{i}<0
\end{cases}\label{eq:SE3PPF_Smooth}
\end{equation}

For simplicity, let $\bar{\delta}=[\bar{\delta}_{1},\ldots,\bar{\delta}_{4}]^{\top}$,
$\underline{\delta}=[\underline{\delta}_{1},\ldots,\underline{\delta}_{4}]^{\top}$,
$\mathcal{E}=\left[\mathcal{E}_{{\rm R}},\mathcal{E}_{{\rm P}}^{\top}\right]^{\top}$
for all $\bar{\delta},\underline{\delta},\mathcal{E}\in\mathbb{R}^{4}$
with $\mathcal{E}_{{\rm R}}=\mathcal{E}_{1}\in\mathbb{R}$ and $\mathcal{E}_{{\rm P}}=[\mathcal{E}_{2},\mathcal{E}_{3},\mathcal{E}_{4}]^{\top}\in\mathbb{R}^{3}$.
The inverse transformation in \eqref{eq:SE3PPF_Smooth} is equivalent
to
\begin{equation}
\begin{aligned}\mathcal{E}_{i}= & \frac{1}{2}\begin{cases}
\text{ln}\frac{\underline{\delta}_{i}+\boldsymbol{e}_{i}/\xi_{i}}{\bar{\delta}_{i}-\boldsymbol{e}_{i}/\xi_{i}}, & \bar{\delta}_{i}\geq\underline{\delta}_{i}\text{ if }\boldsymbol{e}_{i}\geq0\\
\text{ln}\frac{\underline{\delta}_{i}+\boldsymbol{e}_{i}/\xi_{i}}{\bar{\delta}_{i}-\boldsymbol{e}_{i}/\xi_{i}}, & \underline{\delta}_{i}\geq\bar{\delta}_{i}\text{ if }\boldsymbol{e}_{i}<0
\end{cases}\end{aligned}
\label{eq:SE3PPF_trans2}
\end{equation}

\begin{rem}
	\label{rem:SO3PPF_1} Consider $\mathcal{E}_{i}$ in \eqref{eq:SE3PPF_trans2}.
	$\boldsymbol{e}_{i}$ is bounded by $\xi_{i}$, and the prescribed
	performance is achieved if and only if $\mathcal{E}_{i}$ is bounded
	for all $t\geq0$ and $i=1,2,\ldots,4$. 
\end{rem}
\begin{prop}
	\label{Prop:SE3PPF_1}Consider the transformed error in \eqref{eq:SE3PPF_trans2}
	with $\underline{\delta}=\bar{\delta}$, then the following statements
	hold: 
	\begin{enumerate}
		\item[(i)] $\mathcal{E}=\underline{\mathbf{0}}_{4}$ only at $\boldsymbol{e}=\underline{\mathbf{0}}_{4}$,
		and the critical point of $\mathcal{E}$ satisfies $\boldsymbol{e}=\underline{\mathbf{0}}_{4}$. 
		\item[(ii)] The only critical point of $\mathcal{E}$ is $\tilde{\boldsymbol{T}}=\mathbf{I}_{4}$. 
	\end{enumerate}
\end{prop}
\textbf{Proof. }Since $\underline{\delta}=\bar{\delta}$ with the
constraint of $\boldsymbol{e}_{i}\leq\xi_{i}$, it is obvious that
$(\underline{\delta}_{i}+\boldsymbol{e}_{i}/\xi_{i})/(\bar{\delta}_{i}-\boldsymbol{e}_{i}/\xi_{i})=1$
only at $\boldsymbol{e}_{i}=0$. Thus, $\mathcal{E}_{i}\neq0\forall\boldsymbol{e}_{i}\neq0$
and $\mathcal{E}_{i}=0$ only at $\boldsymbol{e}_{i}=0$ which proves
(i). For (ii), from \eqref{eq:SE3PPF_Terr_matrix}, $||\tilde{R}||_{I}=0$
and $\tilde{P}=0$ only at $\tilde{\boldsymbol{T}}=\mathbf{I}_{4}$.
Therefore, the critical point of $\mathcal{E}$ satisfies $||\tilde{R}||_{I}=0$
and $\tilde{P}=0$ which implies that $\tilde{\boldsymbol{T}}=\mathbf{I}_{4}$
and justifies (ii). Define $\mu_{i}:=\mu_{i}\left(\boldsymbol{e}_{i},\xi_{i}\right)$
by
\begin{equation}
\begin{split}\mu_{i} & =\frac{1}{2\xi_{i}}(\frac{1}{\underline{\delta}_{i}+\boldsymbol{e}_{i}/\xi_{i}}+\frac{1}{\bar{\delta}_{i}-\boldsymbol{e}_{i}/\xi_{i}})\end{split}
\label{eq:SE3PPF_mu}
\end{equation}
Let $x=\frac{\dot{\xi}_{1}}{\xi_{1}}$, $X={\rm diag}\left(\frac{\dot{\xi}_{2}}{\xi_{2}},\frac{\dot{\xi}_{3}}{\xi_{3}},\frac{\dot{\xi}_{4}}{\xi_{4}}\right)$,
and $\mathcal{M}={\rm diag}(\mu_{2},\mu_{3},\mu_{4})$ for all $x\in\mathbb{R}$
and $X,\mathcal{M}\in\mathbb{R}^{3\times3}$. Hence, it can be found
that{\small{}
	\begin{equation}
	\dot{\mathcal{E}}=\left[\begin{array}{cc}
	\mu_{1} & \underline{\mathbf{0}}_{3}^{\top}\\
	\underline{\mathbf{0}}_{3} & \mathcal{M}
	\end{array}\right]\left(\left[\begin{array}{c}
	||\dot{\tilde{R}}||_{I}\\
	\dot{\tilde{P}}
	\end{array}\right]-\left[\begin{array}{cc}
	x & \underline{\mathbf{0}}_{3}^{\top}\\
	\underline{\mathbf{0}}_{3} & X
	\end{array}\right]\left[\begin{array}{c}
	||\tilde{R}||_{I}\\
	\tilde{P}
	\end{array}\right]\right)\label{eq:SE3PPF_Trans_dot}
	\end{equation}
}The following section presents a nonlinear pose filter on $\mathbb{SE}\left(3\right)$
with prescribed performance guaranteeing $\mathcal{E}_{i}\in\mathcal{L}_{\infty},\forall t\geq0$.

\section{{\normalsize{}Nonlinear Pose Filter On $\mathbb{SE}\left(3\right)$
		with Prescribed Performance \label{sec:SE3PPF-Filters}}}

This section presents a nonlinear complementary pose filter on $\mathbb{SE}\left(3\right)$
with the error vector in \eqref{eq:SE3PPF_Vec_error} following transient
as well as steady-state measures predefined by the user. Consider
the error in \eqref{eq:SE3PPF_Vec_error}. Define $\boldsymbol{T}_{y}=\left[\begin{array}{cc}
R_{y} & P_{y}\\
\underline{\mathbf{0}}_{3}^{\top} & 1
\end{array}\right]$ as a reconstructed homogeneous transformation matrix of the true
$\boldsymbol{T}$. $R_{y}$ which is corrupted with uncertainty in
measurements is reconstructed using singular value decomposition \cite{markley1988attitude},
or for simplicit visit the appendix in \cite{hashim2018SO3Stochastic}.
$P_{y}$ is reconstructed by 
\[
P_{y}=\frac{1}{\sum_{j=1}^{N_{{\rm L}}}k_{j}^{{\rm L}}}\sum_{j=1}^{N_{{\rm L}}}k_{j}^{{\rm L}}\left({\rm v}_{j}^{\mathcal{I}\left({\rm L}\right)}-R_{y}{\rm v}_{j}^{\mathcal{B}\left({\rm L}\right)}\right)
\]
Consider the following pose filter design
\begin{align}
\dot{\hat{\boldsymbol{T}}}= & \left[\begin{array}{cc}
\hat{R} & \hat{P}\\
\underline{\mathbf{0}}_{3}^{\top} & 1
\end{array}\right]\left[\begin{array}{cc}
\left[\hat{\Omega}\right]_{\times} & \hat{V}\\
\underline{\mathbf{0}}_{3}^{\top} & 0
\end{array}\right]\label{eq:SE3PPF_Test_dot_Ty}\\
\hat{\Omega}= & \Omega_{m}-\hat{b}_{\Omega}-\hat{R}^{\top}W_{\Omega}\label{eq:SE3PPF_Omega_est_Ty}\\
\hat{V}= & V_{m}-\hat{b}_{V}+W_{V}\label{eq:SE3PPF_Vest_Ty}
\end{align}
\begin{align}
W= & \left[\begin{array}{cc}
2\frac{k_{w}\mu_{1}\mathcal{E}_{{\rm R}}-x/4}{1-||\tilde{R}||_{I}} & \mathbf{0}_{3\times3}\\
\mathbf{0}_{3\times3} & k_{w}\hat{R}^{\top}
\end{array}\right]\left[\begin{array}{c}
\mathbf{vex}(\boldsymbol{\mathcal{P}}_{a}(\tilde{R}))\\
\mathcal{M}\mathcal{E}_{{\rm P}}
\end{array}\right]\nonumber \\
& +\left[\begin{array}{c}
\mathbf{0}_{3\times1}\\
\left[\tilde{P}-\hat{P}\right]_{\times}W_{\Omega}-X\tilde{P}
\end{array}\right]\label{eq:SE3PPF_Wcorr_Ty}
\end{align}
\begin{equation}
\dot{\hat{b}}=\gamma\left[\begin{array}{cc}
\frac{1}{2}\mu_{1}\mathcal{E}_{{\rm R}}\hat{R}^{\top} & \hat{R}^{\top}\left[\tilde{P}-\hat{P}\right]_{\times}\\
\mathbf{0}_{3\times3} & \hat{R}^{\top}
\end{array}\right]\left[\begin{array}{c}
\mathbf{vex}(\boldsymbol{\mathcal{P}}_{a}(\tilde{R}))\\
\mathcal{M}\mathcal{E}_{{\rm P}}
\end{array}\right]\label{eq:SE3PPF_best_dot_Ty}
\end{equation}
with $\tilde{R}=\hat{R}R_{y}^{\top}$, $\tilde{P}=\hat{P}-\tilde{R}P_{y}$,
$\mathcal{E}_{{\rm R}}$, $\mathcal{E}_{{\rm P}}$, $\mu_{1}$ and
$\mathcal{M}$ being defined in \eqref{eq:SE3PPF_mu} and \eqref{eq:SE3PPF_Trans_dot},
$k_{w}$ and $\gamma$ being positive constants, $W=\left[W_{\Omega}^{\top},W_{V}^{\top}\right]^{\top}$
and $\hat{b}=\left[\hat{b}_{\Omega}^{\top},\hat{b}_{V}^{\top}\right]^{\top}$
being the correction factor and the estimate of $b$, respectively.
Define the error between the true and the estimated bias by
\begin{align}
\tilde{b} & =b-\hat{b}\label{eq:SE3PPF_b1_tilde}
\end{align}
with $\tilde{b}=\left[\tilde{b}_{\Omega}^{\top},\tilde{b}_{V}^{\top}\right]^{\top}\in\mathbb{R}^{6}$
being the group error bias vector.  
\begin{thm}
	\textbf{\label{thm:SE3PPF_1} }Consider the pose kinematics in \eqref{eq:SE3PPF_T_Dynamics}
	and the group of noise-free velocity measurements in \eqref{eq:SE3PPF_Angular}
	where $\mathcal{Y}_{m}=\Omega+b$, in addition to other vector measurements
	given in \eqref{eq:SE3STCH_Vector_norm} and \eqref{eq:SE3STCH_Vec_Landmark}
	coupled with the filter in \eqref{eq:SE3PPF_Test_dot_Ty}, \eqref{eq:SE3PPF_Omega_est_Ty},
	\eqref{eq:SE3PPF_Vest_Ty}, \eqref{eq:SE3PPF_Wcorr_Ty}, and \eqref{eq:SE3PPF_best_dot_Ty}.
	Let Assumption \ref{Assum:SE3STCH_1} hold. Define $\mathcal{U}\subseteq\mathbb{SE}\left(3\right)\times\mathbb{R}^{6}$
	by{\small{} $\mathcal{U}:=\left\{ \left.(\tilde{\boldsymbol{T}}\left(0\right),\tilde{b}\left(0\right))\right|{\rm Tr}\{\tilde{R}\left(0\right)\}=-1,\tilde{P}\left(0\right)=\underline{\mathbf{0}}_{3},\tilde{b}\left(0\right)=\underline{\mathbf{0}}_{6}\right\} $}.
	For $\tilde{R}\left(0\right)\notin\mathcal{U}$, $\tilde{P}\in \mathbb{R}^{3}$ and $\mathcal{E}\left(0\right)\in\mathcal{L}_{\infty}$, all the closed loop signals are uniformly ultimately bounded,
	$\lim_{t\rightarrow\infty}\mathcal{E}\left(t\right)=0$ and $\tilde{\boldsymbol{T}}$
	asymptotically approaches $\mathbf{I}_{4}$.  
\end{thm}
\textbf{Proof. }Let the error in the homogeneous transformation matrix
be as in \eqref{eq:SE3PPF_Terr_matrix}. From \eqref{eq:SE3PPF_T_Dynamics}
and \eqref{eq:SE3PPF_Test_dot_Ty} the error in attitude dynamics
is
\begin{align}
\dot{\tilde{R}} & =\left[\hat{R}\tilde{b}_{\Omega}-W_{\Omega}\right]_{\times}\tilde{R}\label{eq:SE3PPF_Rtilde_dot}
\end{align}
In view of \eqref{eq:SE3PPF_T_Dynamics} and \eqref{eq:SE3PPF_NormR_dynam},
the error dynamics in \eqref{eq:SE3PPF_Rtilde_dot} can be expressed
in terms of the normalized Euclidean distance
\begin{align}
\frac{d}{dt}||\tilde{R}||_{I} & =\frac{1}{2}\mathbf{vex}(\boldsymbol{\mathcal{P}}_{a}(\tilde{R}))^{\top}(\hat{R}\tilde{b}_{\Omega}-W_{\Omega})\label{eq:SE3PPF_NormRtilde_dot}
\end{align}
where ${\rm Tr}\left\{ \tilde{R}\left[W_{\Omega}\right]_{\times}\right\} =-2\mathbf{vex}(\boldsymbol{\mathcal{P}}_{a}(\tilde{R}))^{\top}W_{\Omega}$
as in \eqref{eq:SO3PPF_Identity7}. The derivative of $\tilde{P}$
can be found to be
\begin{align}
\dot{\tilde{P}} & =\hat{R}(\tilde{b}_{V}-W_{V})+\left[\hat{P}-\tilde{P}\right]_{\times}(\hat{R}\tilde{b}_{\Omega}-W_{\Omega})\label{eq:SE3PPF_Ptilde_dot}
\end{align}
with $\left[\hat{R}\tilde{b}_{\Omega}\right]_{\times}\hat{P}=-\left[\hat{P}\right]_{\times}\hat{R}\tilde{b}_{\Omega}$.
The dynamics of the error vector in \eqref{eq:SE3PPF_Vec_error} become{\small{}
	\begin{align}
	\left[\begin{array}{c}
	||\dot{\tilde{R}}||_{I}\\
	\dot{\tilde{P}}
	\end{array}\right] & =\left[\begin{array}{cc}
	\frac{1}{2}\mathbf{vex}(\boldsymbol{\mathcal{P}}_{a}(\tilde{R}))^{\top} & \underline{\mathbf{0}}_{3}^{\top}\\
	\left[\hat{P}-\tilde{P}\right]_{\times} & \hat{R}
	\end{array}\right]\left[\begin{array}{c}
	\hat{R}\tilde{b}_{\Omega}-W_{\Omega}\\
	\tilde{b}_{V}-W_{V}
	\end{array}\right]\label{eq:SE3PPF_VEC_tilde_dot}
	\end{align}
}Consider the following candidate Lyapunov function
\begin{align}
V(\mathcal{E},\tilde{b}) & =\frac{1}{2}||\mathcal{E}||^{2}+\frac{1}{2\gamma}||\tilde{b}||^{2}\label{eq:SE3PPF_V_Ry}
\end{align}
Differentiating $V:=V(\mathcal{E},\tilde{b})$ in \eqref{eq:SE3PPF_V_Ry},
and considering $||\tilde{R}||_{I}=\frac{1}{4}\frac{\left\Vert \mathbf{vex}(\boldsymbol{\mathcal{P}}_{a}(\tilde{R}))\right\Vert ^{2}}{1-||\tilde{R}||_{I}}$
as defined in \eqref{eq:SE3PPF_lemm1_1} with direct substitution
of $\dot{\hat{b}}$ and $W$ in \eqref{eq:SE3PPF_best_dot_Ty} and
\eqref{eq:SE3PPF_Wcorr_Ty}, respectively, one obtains
\begin{align}
\dot{V} & =-k_{w}\mathcal{E}^{\top}\left[\begin{array}{cc}
||\tilde{R}||_{I} & \underline{\mathbf{0}}_{3}^{\top}\\
\underline{\mathbf{0}}_{3} & \mathbf{I}_{3}
\end{array}\right]\left[\begin{array}{cc}
\mu_{1} & \underline{\mathbf{0}}_{3}^{\top}\\
\underline{\mathbf{0}}_{3} & \mathcal{M}
\end{array}\right]^{2}\mathcal{E}\label{eq:SE3PPF_Vdot_Ry_Final}
\end{align}
The result in \eqref{eq:SE3PPF_Vdot_Ry_Final} indicates that $V\left(t\right)\leq V\left(0\right),\forall t\geq0$
and $\tilde{R}\left(0\right)\notin\mathcal{U}$. Consequently, $\tilde{b}$
and $\mathcal{E}$ remain bounded for all $t\geq0$. Thus, $\tilde{P}$,
$||\tilde{R}||_{I}$ and $\mathbf{vex}(\boldsymbol{\mathcal{P}}_{a}(\tilde{R}))$
are bounded which in turn implies the boundedness of $\dot{\tilde{P}}$,
$||\dot{\tilde{R}}||_{I}$, $\dot{\mathcal{E}}_{{\rm R}}$ and $\dot{\mathcal{E}}_{{\rm P}}$.
One can find that $\dot{\mu}_{i}$ is bounded as well which means
that $\ddot{V}$ is bounded for all $t\geq0$. Therefore, $\dot{V}$
is uniformly continuous and, consistent with Barbalat Lemma, $\dot{V}\rightarrow0$
as $t\rightarrow\infty$ indicating that $\mathcal{E}_{i}\rightarrow0$
and $\boldsymbol{e}_{i}\rightarrow0$ for all $i=1,2,\ldots,4$. According to property (ii)
of Proposition \ref{Prop:SE3PPF_1}, $\mathcal{E}\rightarrow0$ implies
that $\tilde{\boldsymbol{T}}$ asymptotically approaches $\mathbf{I}_{4}$
which completes the proof.

\section{Simulations \label{sec:SO3PPF_Simulations}}

Let the dynamics of $\boldsymbol{T}$ be as described in \eqref{eq:SE3PPF_T_Dynamics}.
Define the true angular and translational velocities by{\small{}
	\begin{align*}
	\Omega & =0.8\left[0.6{\rm sin}\left(0.4t\right),{\rm cos}\left(0.6t\right),0.7{\rm sin}\left(0.3t+\frac{\pi}{5}\right)\right]^{\top}\left({\rm rad/sec}\right)\\
	V & =0.3\left[0.4{\rm cos}\left(0.5t\right),{\rm sin}\left(0.2t\right),0.2{\rm sin}\left(0.4t+\frac{\pi}{3}\right)\right]^{\top}\left({\rm m/sec}\right)
	\end{align*}
}Let $\Omega_{m}=\Omega+b_{\Omega}+\omega_{\Omega}$ and $V_{m}=V+b_{V}+\omega_{V}$,
with $b_{\Omega}=0.1\left[1,-1,1\right]^{\top}$ and $b_{V}=0.1\left[2,5,1\right]^{\top}$.
$\omega_{\Omega}$ and $\omega_{V}$ represent random noise with zero
mean and standard deviation (STD) equal to $0.16\left({\rm rad/sec}\right)$
and $0.25\left({\rm m/sec}\right)$, respectively. Assume $N_{{\rm L}}=1$
and $N_{{\rm R}}=2$ with ${\rm v}_{1}^{\mathcal{I}\left({\rm L}\right)}=\left[\frac{1}{2},\sqrt{2},1\right]^{\top}$,
${\rm v}_{1}^{\mathcal{I}\left({\rm R}\right)}=\frac{1}{\sqrt{3}}\left[1,-1,1\right]^{\top}$
and ${\rm v}_{2}^{\mathcal{I}\left({\rm R}\right)}=\left[0,0,1\right]^{\top}$.
Let ${\rm v}_{1}^{\mathcal{B}\left({\rm L}\right)}=R^{\top}\left({\rm v}_{1}^{\mathcal{I}\left({\rm L}\right)}-P\right)+{\rm b}_{1}^{\mathcal{B}\left({\rm L}\right)}+\omega_{1}^{\mathcal{B}\left({\rm L}\right)}$
and ${\rm v}_{j}^{\mathcal{B}\left({\rm R}\right)}=R^{\top}{\rm v}_{j}^{\mathcal{I}\left({\rm R}\right)}+{\rm b}_{j}^{\mathcal{B}\left({\rm R}\right)}+\omega_{j}^{\mathcal{B}\left({\rm R}\right)}$
for $j=1,2$ with ${\rm b}_{1}^{\mathcal{B}\left({\rm L}\right)}=0.1\left[0.3,0.2,-0.2\right]^{\top}$,
${\rm b}_{1}^{\mathcal{B}\left({\rm R}\right)}=0.1\left[-1,1,0.5\right]^{\top}$
and ${\rm b}_{2}^{\mathcal{B}\left({\rm R}\right)}=0.1\left[0,0,1\right]^{\top}$.
Additionally, $\omega_{1}^{\mathcal{B}\left({\rm L}\right)}$, $\omega_{1}^{\mathcal{B}\left({\rm R}\right)}$
and $\omega_{2}^{\mathcal{B}\left({\rm R}\right)}$ are Gaussian noise
vectors with zero mean and ${\rm STD}=0.3$, ${\rm STD}=0.1$, and
${\rm STD}=0.1$, respectively. In order to satisfy Assumption \ref{Assum:SE3STCH_1},
the third vector is obtained using ${\rm v}_{3}^{\mathcal{I}\left({\rm R}\right)}={\rm v}_{1}^{\mathcal{I}\left({\rm R}\right)}\times{\rm v}_{2}^{\mathcal{I}\left({\rm R}\right)}$
and ${\rm v}_{3}^{\mathcal{B}\left({\rm R}\right)}={\rm v}_{1}^{\mathcal{B}\left({\rm R}\right)}\times{\rm v}_{2}^{\mathcal{B}\left({\rm R}\right)}$.
Next, ${\rm v}_{j}^{\mathcal{B}\left({\rm R}\right)}$ and ${\rm v}_{j}^{\mathcal{I}\left({\rm R}\right)}$
are normalized to $\upsilon_{j}^{\mathcal{B}\left({\rm R}\right)}$
and $\upsilon_{j}^{\mathcal{I}\left({\rm R}\right)}$, respectively,
for $j=1,2,3$ as given in \eqref{eq:SE3STCH_Vector_norm}. $R_{y}$
is obtained by SVD (visit the appendix in \cite{hashim2018SO3Stochastic}
or \cite{hashim2018SE3Stochastic}) with $\tilde{R}=\hat{R}R_{y}^{\top}$.
The initialization of the true and the estimated pose is given by{\small{}
	\[
	\boldsymbol{T}\left(0\right)=\mathbf{I}_{4},\hspace{1em}\hat{\boldsymbol{T}}\left(0\right)=\left[\begin{array}{cccc}
	-0.8816 & 0.2386 & 0.4074 & -4\\
	0.4498 & 0.1625 & 0.8782 & 5\\
	0.1433 & 0.9574 & -0.2505 & 3\\
	0 & 0 & 0 & 1
	\end{array}\right]
	\]
}The design parameters of the proposed filters are chosen as $\gamma=1$,
$k_{w}=6$, $\bar{\delta}=\underline{\delta}=\left[1.3,5,6,4\right]^{\top}$,
$\xi^{0}=\left[1.3,5,6,4\right]^{\top}$, $\xi^{\infty}=\left[0.07,0.3,0.3,0.3\right]^{\top}$,
and $\ell=\left[4,4,4,4\right]^{\top}$. The initial bias estimate
is $\hat{b}\left(0\right)=\underline{\mathbf{0}}_{6}$.

\begin{figure}[h]
	\centering{}\includegraphics[scale=0.27]{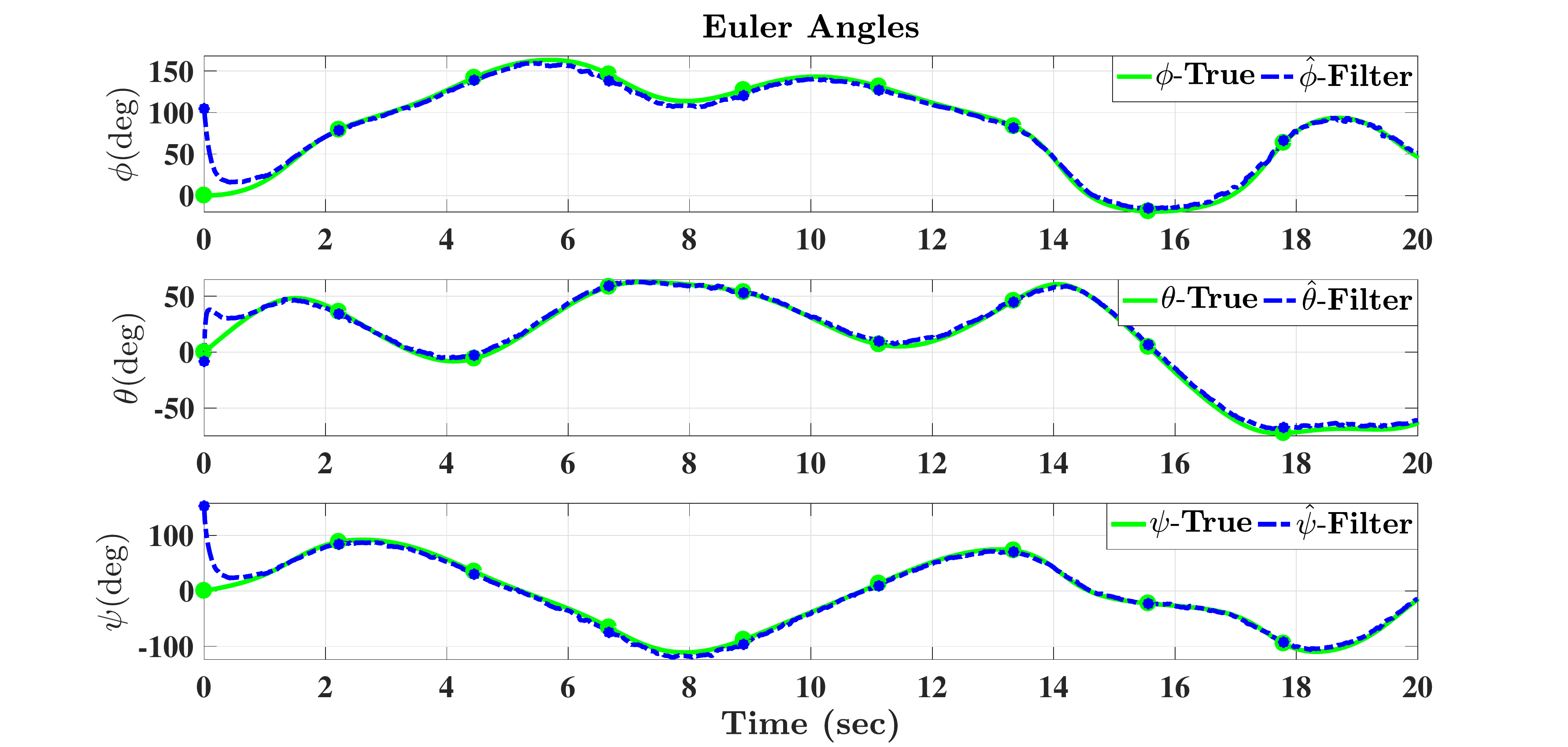}\caption{True and estimated Euler angles of the rigid-body.}
	\label{fig:SE3PPF_Simulation4}  
\end{figure}

\begin{figure}[h]
	\centering{}\includegraphics[scale=0.27]{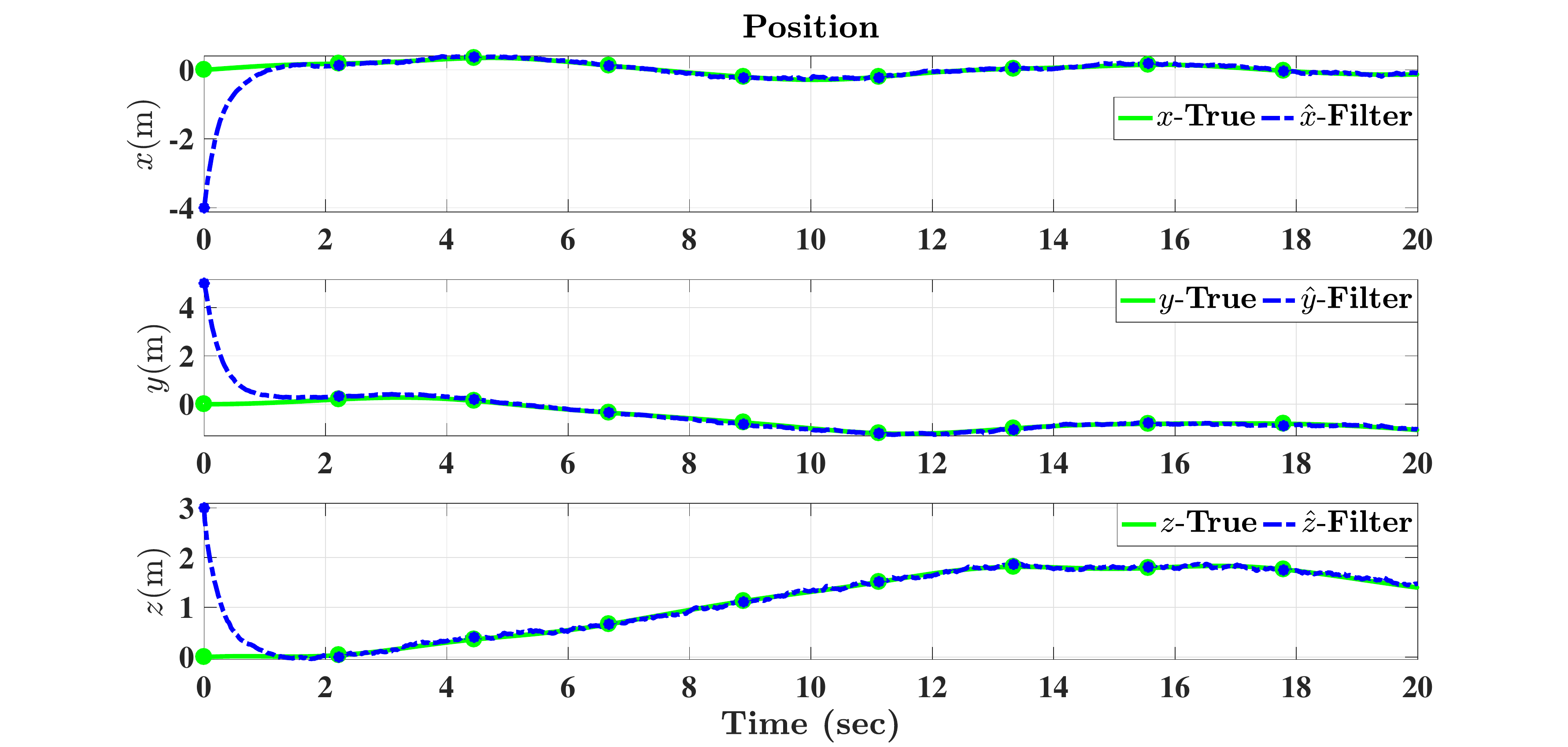}\caption{True and estimated position of rigid-body in 3D space.}
	\label{fig:SE3PPF_Simulation5}  
\end{figure}

\begin{figure*}[h]
	\centering{}\includegraphics[scale=0.52]{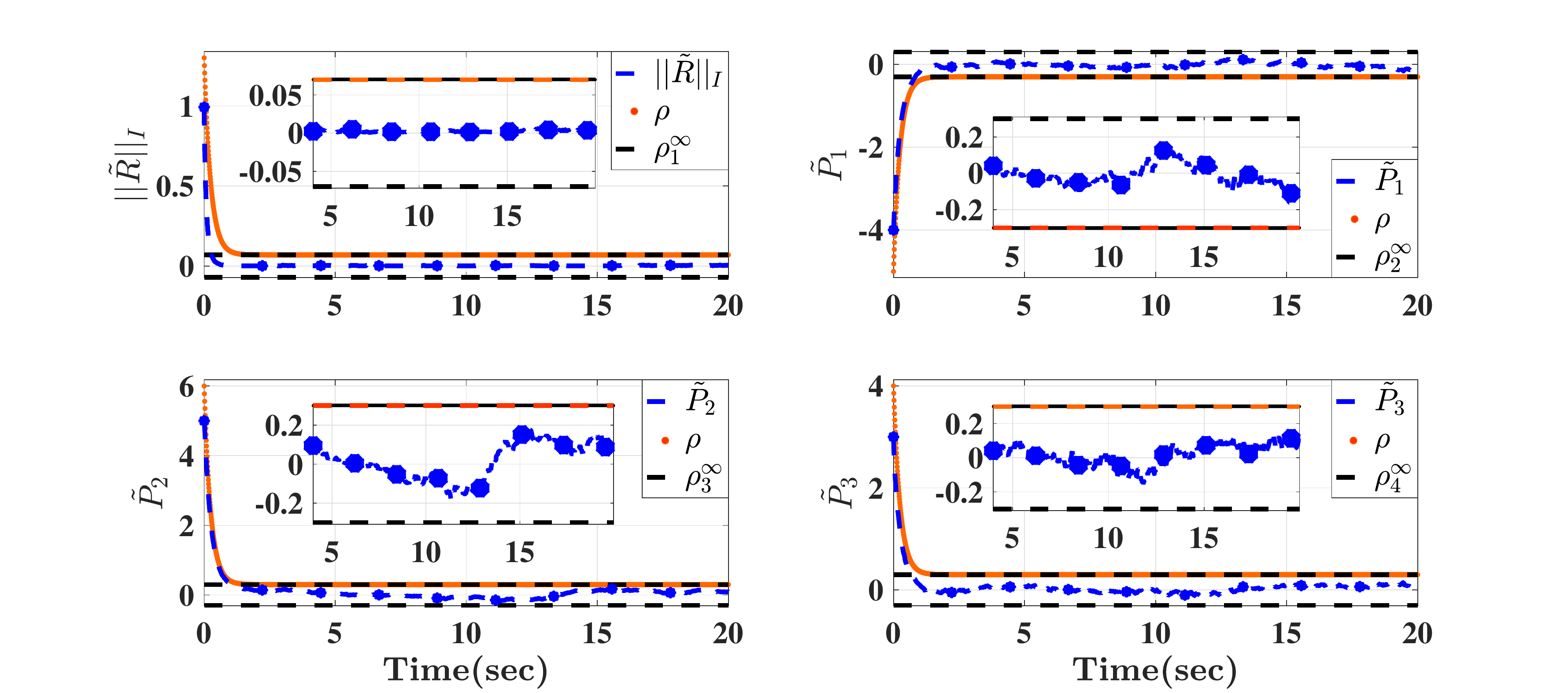}\caption{Systematic convergence of the error trajectories within the dynamic
		decreasing boundaries.}
	\label{fig:SE3PPF_Simulation6}  
\end{figure*}

Fig \ref{fig:SE3PPF_Simulation4} and \ref{fig:SE3PPF_Simulation5}
show impressive tracking performance with fast convergence of the
Euler angles $\left(\phi,\theta,\psi\right)$ and $xyz$-coordinates
in 3D space, respectively. Fig. \ref{fig:SE3PPF_Simulation6} illustrates
the systematic and smooth convergence of the error vector $\boldsymbol{e}$
demonstrating that $||\tilde{R}||_{I}$ starts very close to the unstable
equilibria ($+1$) while $\tilde{P}_{1}$, $\tilde{P}_{2}$, and $\tilde{P}_{3}$
start with large error within the predefined large set and attenuate
systematically to the predefined small set.

\section{Conclusion\label{sec:SO3PPF_Conclusion}}

A nonlinear pose filter with predefined characteristics has been introduced.
The filter is evolved directly on $\mathbb{SE}\left(3\right)$. The
pose error has been formulated in terms of position error and normalized
Euclidean distance error. The error vector has been constrained to
follow the predefined dynamically decreasing boundaries such that
the transient performance does not exceed the dynamically decreasing
function and the error is regulated to the origin asymptotically from
almost any initial condition. %
{} Simulation results showed robustness of the proposed filter against
high level of uncertainties in the measurements and large initialization
error.

\section*{Acknowledgment}

The authors would like to thank \textbf{Maria Shaposhnikova} for proofreading
the article.

\section*{Appendix A \label{sec:SO3STCH_EXPL_AppendixA}}
\begin{center}
	\textbf{Proof of Lemma \ref{Lemm:SE3PPF_1}}
	\par\end{center}

Define $R\in\mathbb{SO}\left(3\right)$ as an attitude of a rigid-body
in 3D space. The attitude can be represented in terms of Rodriguez
parameters vector $\rho\in\mathbb{R}^{3}$ while the mapping from
$\rho$ to $\mathbb{SO}\left(3\right)$ is defined by $\mathcal{R}_{\rho}:\mathbb{R}^{3}\rightarrow\mathbb{SO}\left(3\right)$
\cite{shuster1993survey,hashim2018SO3Stochastic} 
\begin{align}
\mathcal{R}_{\rho}\left(\rho\right)= & \frac{1}{1+||\rho||^{2}}\left(\left(1-||\rho||^{2}\right)\mathbf{I}_{3}+2\rho\rho^{\top}+2\left[\rho\right]_{\times}\right)\label{eq:SE3PPF_SO3_Rodr}
\end{align}
substituting \eqref{eq:SE3PPF_SO3_Rodr} into \eqref{eq:SE3STCH_Ecul_Dist}
one has 
\begin{equation}
||R||_{I}=\frac{||\rho||^{2}}{1+||\rho||^{2}}\label{eq:SE3PPF_TR2}
\end{equation}
The anti-symmetric projection operator in \eqref{eq:SE3PPF_SO3_Rodr}
is
\begin{align}
\boldsymbol{\mathcal{P}}_{a}\left(R\right)=\frac{1}{2}\left(\mathcal{R}_{\rho}-\mathcal{R}_{\rho}^{\top}\right)= & 2\frac{1}{1+||\rho||^{2}}\left[\rho\right]_{\times}\label{eq:SE3PPF_Pa}
\end{align}
As such, the vex operator in \eqref{eq:SE3PPF_Pa} is
\begin{equation}
\mathbf{vex}\left(\boldsymbol{\mathcal{P}}_{a}\left(R\right)\right)=2\frac{\rho}{1+||\rho||^{2}}\label{eq:SE3PPF_VEX_Pa}
\end{equation}
From \eqref{eq:SE3PPF_TR2} one finds
\begin{equation}
\left(1-||R||_{I}\right)||R||_{I}=\frac{||\rho||^{2}}{\left(1+||\rho||^{2}\right)^{2}}\label{eq:SE3PPF_append1}
\end{equation}
and from \eqref{eq:SE3PPF_VEX_Pa} one has
\begin{equation}
||\mathbf{vex}\left(\boldsymbol{\mathcal{P}}_{a}\left(R\right)\right)||^{2}=4\frac{||\rho||^{2}}{\left(1+||\rho||^{2}\right)^{2}}\label{eq:SE3PPF_append2}
\end{equation}
Thus, \eqref{eq:SE3PPF_append1} and \eqref{eq:SE3PPF_append2} justify
\eqref{eq:SE3PPF_lemm1_1} in Lemma \ref{Lemm:SE3PPF_1}. 

\bibliographystyle{IEEEtran}
\bibliography{bib_PPF_SE3}

\end{document}